\documentclass[12pt]{amsart}
\usepackage{latexsym}
\usepackage{amsfonts}
\usepackage{amsthm}
\usepackage{amsmath}
\usepackage{amssymb}

\def\<{{\langle}}
\def\>{{\rangle}}

\def\eps{\varepsilon}

\def\note#1{{}}

\def\note#1{}

\def\beq{\begin{equation}}
\def\eeq{\end{equation}}

\newcounter{zlist}

\def\Label#1{\label{#1}\ifmmode\llap{[#1] }\else
\marginpar{\smash{\hbox{\tiny [#1]}}}\fi}
\def\Label{\label}

\newtheorem{proposition}{Proposition}[section]

\newtheorem{theorem}[proposition]{Theorem}

\theoremstyle{definition}
\newtheorem{definition}[proposition]{Definition}

\theoremstyle{remark}
\newtheorem{remark}[proposition]{Remark}

\newcounter{c}

\newcommand{\etyk}[1]{\vspace{-7.4mm}$$\begin{equation}\Label{#1}
\addtocounter{c}{1}}
\renewcommand{\]}{\ifnum \value{c}=1 $$\else \end{equation}\fi}
\begin{document}

\title{Memories with Solomon Marcus (III)}
\author{Florin Felix NICHITA}
\address{Institute of Mathematics of the Romanian Academy, 
P.O. Box 1-764, RO-70700, Bucharest, Romania}
\date{}
\begin{abstract} { Dedicated
to Solomon Marcus, the current review
 adds interesting details
to our previous articles. Some new results
are also included.
 I was interested in
the work of Solomon Marcus
in Mathematical Linguistics as a high-school student.
Later,
I had the opportunity  to discuss with him
about many topics. He was a polymath.
We wrote a paper together, and
I refereed an  editorial paper
about his work in 2021.
Samples of  (possible)
discussions are presented: some topology conjectures,
a self-dual theorem in geometry,
results about Boolean algebras,  
 a B--ring Euler formula, Yang-Baxter maps
 and a discussion on sequences and series.
A short appendix on poetry is also included.
}
\end{abstract}
\maketitle

\section{Introduction} 

After writing { \em Memories with Solomon Marcus (I)} and { \em (II)},
I felt that we need an unified version of them with further developments
included.
Even from  high-school, I interested in the
work of Solomon Marcus (see \cite{11}) in Mathematical Linguistics.
As I had participated in some
communication sessions for high-school students,  
I heard about the entropy of a translated poetry (see \cite{22})
from and older colleague. 
My colleague and me, won good prizes. From those manuscripts,
I publishing an exercise
in the
Romanian Gazette (see \cite{Memo}). 
Several years later, I wrote an article in the
 Romanian Mathematical Gazette on that exercise (see \cite{A}). 
Some math teachers from Bucharest wrote
a book including results from that article. 
I had the opportunity to meet Solomon Marcus  and discuss with him 
as a student and as a researcher at 
 the
Institute of Mathematics of the Romanian Academy
(see \cite{22}).
As a result of our ``unexpected'' meetings, we wrote a paper 
together (see \cite{SM}). We got some of these ideas after a long
discussion on Euclidean geometry. It might be the case that Solomon Marcus
was actually impressed by the paper \cite{Euler}, and he also wrote
on transcendental numbers several papers during those years.
I will refer to some remarks of \cite{SM} in the next section. 
Solomon Marcus also invited my friends and me to many events. I spent 
many hours
in the libraries studying his works in Boolean algebras, informatics,
literature and philosophy, trying to give a replica
to his legacy (see \cite{F,G,H}). When I met one of the most
important Romanian writers, Mircea Cartarescu, I asked him: 
``What do you think
about the connection between Mathematics and Literature?''
He only said to me: ``Solomon Marcus''...
In 2021,
I was editor for an article  about Solomon Marcus' work (\cite{editor}).

\bigskip

\section{Overture}

In the present paper, we recall  some  previous results 
 (from
\cite{Memo, mMemo}), adding some  interesting remarks and new theorems.
This part of the paper could be considered a continuation
of our Axioms paper (\cite{SM}).
The next two sections contain
conjectures in topology and unification constructions in geometry.
Section 5 is on Boolean Algebra Theory. We try to cast
the main properties of the Euler formula in the Boolean algebra setting.
Even in this rigid environment we are able to obtain
interesting results. We will need to relax the restricting
conditions imposed by the Boolean algebras category in order
to obtain finer results. So,
we  present the concept of (weak) B-rings,
an unifying structure for Boolean algebras and rings.
The main  B-rings for our theory will arise from Boolean algebras.
We  will find 
an Euler formula for B-rings. 
There will be a section on (arithmetic / geometric) series.
Started as a high-school problem it lead to a sophisticated
approach to unify the arithmetic and geometric series.
An  appendix on mathematics and poetry will also be included.

\bigskip

\section{Topology}

After a series of conferences 
on Intersection Cohomology,
I remembered some conjectures from
the paper \cite{SM}, and I improvised  a beautiful problem,
which was well-received. 
I will refer  to enhanced versions of our
conjectures below. 
Solomon Marcus and me formulated them after a long
discussion on Euclidean geometry in a late December, just
before Christmas.

\begin{remark}
For a Jordan curve in the Euclidean plane, 
we consider the maximum diameter  (D)
and the smallest diameter passing through the 
corresponding center of mass ($d$).

(i) If $L$ is the length of the given curve then:
$ \ \  \frac{L}{D} \leq \pi \leq \frac{L}{d} \ \ .$

(ii) Moreover, the first inequality becomes equality
 if and only if the second inequality becomes equality
if and only if the given curve is a circle.

(iii) If the area of the domain inside the
given curve is $A$, then $  \ \ \ d \ D  \ > A \ $.

(iv) The equation $ x^2 - \frac{L}{2} x + A =0 $ is not completely solved; for example,
if the given curve is an ellipse, solving this equation is an unsolved problem. 
\end{remark}


The last Congress of Romanian Mathematicians 
attended by Solomon Marcus was the 
Eighth Congress of Romanian Mathematicians, 
 Iasi, Romania, 26 June - 1 July, 2015.
I remembers our happy discussions at 
this congress.
I  addressed
 an interesting
question to Alex Suciu, and it turned out that it  could lead
to a paper in collaboration with Barbu Berceanu.

\bigskip

\section{Geometry}


In general a theorem in Euclidean geometry is different from its dual.
This happens because in the dual theorem we replace the points with lines
and the lines become points. In this section we will find a counterexample
for this general rule.
As
Solomon Marcus was passionate about the Euclidean geometry,
the following theorem could have been included in our meeting discussions.

{Let us fix the terminology first.
Let $ \mathcal{A} $ and  $ \mathcal{B} $ be two triangles with 
vertices $ \ A_1 , \ 
A_2  $ and $ \ A_3 $, and, respectively
with vertices $ \ B_1 , \ 
B_2  $ and $ \ B_3 $.
We will write $ \ \mathcal{A} \ni\in \mathcal{B} $ if and only if
$ \ A_1 \in B_2 B_3 $ and $ \ B_1 \in A_2 A_3 $.}

\begin{theorem} (\cite{Memo})
 Let
$ \ \mathcal{O} \ni \in \mathcal{X} $,
$ \ \mathcal{O} \ni \in \mathcal{Y} $ and
$ \ \mathcal{O} \ni \in \mathcal{Z} $.
 Let
$ \ \mathcal{X} \ni \in \mathcal{A} $,
$ \ \mathcal{X} \ni \in \mathcal{A'} $,
$ \ \ \ \mathcal{Y} \ni \in \mathcal{B} $,
$ \ \mathcal{Y} \ni \in \mathcal{B '} $,
$ \ \ \ \mathcal{Z} \ni \in \mathcal{C}$  
and
$ \ \mathcal{Z} \ni \in \mathcal{C '} $.
 Let
$ \ \mathcal{A} \ni \in \mathcal{M} $,
$ \ \mathcal{B} \ni \in \mathcal{M} $,
$ \ \mathcal{A'} \ni \in \mathcal{M'} $,
$ \ \mathcal{B'} \ni \in \mathcal{M'} $,
$ \ \mathcal{A} \ni \in \mathcal{N} $,
$ \ \mathcal{C} \ni \in \mathcal{N} $,
$ \ \mathcal{A'} \ni \in \mathcal{N'} $,
$ \ \mathcal{C'} \ni \in \mathcal{N'} $,
$ \ \mathcal{B} \ni \in \mathcal{S} $,
$ \ \mathcal{C} \ni \in \mathcal{S} $,
$ \ \mathcal{B'} \ni \in \mathcal{S'} $ and
$ \ \mathcal{C'} \ni \in \mathcal{S'} $.
 Let
$ \ \mathcal{M} \ni \in \mathcal{P} $,
$ \ \mathcal{M'} \ni \in \mathcal{P} $,
$ \ \  \ \mathcal{N} \ni \in \mathcal{Q} $,
$ \ \mathcal{N'} \ni \in \mathcal{Q} $,
$ \ \ \ \ \mathcal{S} \ni \in \mathcal{R}$  
and
$ \ \mathcal{S'} \ni \in \mathcal{R} $.$ \ \  \ $Then, there exists 
$ \ \mathcal{O'}$ such that 
$ \ \ \  \ \mathcal{P} \ni \in \mathcal{O'} $,
$ \ \mathcal{Q} \ni \in \mathcal{O'} $
and
$ \ \mathcal{R} \ni \in \mathcal{O'} $.

\end{theorem}

\begin{remark}
 This is a self-dual theorem, because the dual of a triangle is
also a triangle.
\end{remark}

\bigskip

\section{Boolean Algebras}

\subsection{A toy model for Euler formula in Boolean Algebras}

Solomon Marcus was very keen on the Euler identity. We had many meetings
in which we have talked about on it.
I am sure he would enjoy the following observations, which try to cast
some properties of the Euler formula in the Boolean algebra setting.

\bigskip

Let $ \ \mathcal{A} \ = \ ( \ A, \vee, \wedge, \ 0, \ 1, \ \hat \ ) \ $ 
be a Boolean algebra.
We consider a new Boolean algebra 
$ \ \mathcal{B} \ =  \ \mathcal{A} \ \times  \  \mathcal{A} \ \times 
 \ \mathcal{A} \ $, 
 with the pointwise operations.
For $ X = (a , \ b , \ c) \ \in \  \mathcal{B} $, 
we define three  functions,
$ \   exp_{B} \ , \ \cos_{B} \ , \ \sin_{B} \ : \ \mathcal{B} \ \rightarrow
\ \mathcal{B} \ , $ by the following formulas:
$$ \ \ \  e^X = exp_{B} ((a , \ b , \ c)) = ( \hat{a},  \hat{b}, \ 1) \ \in \  \mathcal{B} \ , $$

$    \cos(X) = \cos_{B} ((a ,  b ,  c)) = ( 1,  \hat{b}, \ \hat{c}  )  \in   \mathcal{B}  \  $
and 
$  \ \ \sin(X) = \sin_{B} ((a ,  b ,  c)) = ( 0,  0, {c})  \in  \mathcal{B} \ . $

\bigskip

The following properties  hold for the above functions
(see \cite{Memo}):
\begin{equation}\label{expoC}
  e^{X \vee Y}   
=
e^{X} \wedge e^{ Y} \ ; \ \ \ \ 
 \cos (X \vee Y)  =
\cos (X) \wedge \cos ( Y) \ ; 
\end{equation}
\begin{equation}\label{expoCS}
 \sin (X \vee Y) = [ \sin (X) \wedge \cos ( Y) ]  \vee [ \cos (X) \wedge \sin ( Y) ] 
\vee [\sin (X) \wedge \sin ( Y)] \ .
\end{equation}

The following formulas lead to some kind of fundamental trigonometric formula:

$ \sin^{2 \wedge} (X) \vee \cos^{2 \wedge} (X) \vee 
( \cos \circ \sin) (X) = $

$ = [ \sin (X) \wedge \sin (X) ] \vee [ \cos (X) \wedge \cos (X) ] \vee 
( \cos \circ \sin) (X) = $

$ =
\sin (X) \vee  \cos (X)  \vee 
( \cos ( \sin (X)) = $
$ =
( 0,  0, {c}) \vee  ( 1,  \hat{b}, \ \hat{c}  )   \vee 
 \cos  ( 0,  0, {c}) = 
\mathbf{1}  \ . $

Alternatively, the following formula could be considered the 
 fundamental trigonometric formula for this setting:
$$ \sin (X)  \vee 
( \cos \circ \sin) (X) = \mathbf{1}  \ . $$

\begin{theorem} (\cite{Memo})
Let $ J = (0,1,1) \in \  \mathcal{B} $.
Then, the following Euler type formula holds:
\begin{equation}\label{BAE}
 e^{J \wedge X}  = \ \cos (X) \ \vee \ [ \ J \wedge \sin (X) \ ] .
\end{equation}
\end{theorem}

{ \bf Proof.}
The proof of the formula (\ref{BAE}) is direct. Let
$ e^{J \wedge X} = e^{(0,\ b,\ c)} = (1,\ \hat{b}, \ 1) \ $. On the other hand
$ \ \ \cos X \vee [ (0, 1, 1) \wedge \sin X] = 
(1, \hat{b}, \hat{c}) \vee (0, 0, c) = (1, \hat{b}, 1)$.
\qed

\begin{remark} 
Let us denote the element $ \ (1, 0, 0) \in   \mathcal{B} $ by $ \Pi $.
This notation is motivated by the fact that 
$ \ \ \ \ \ \cos (\Pi)  = \mathbf{1} \ 
  \ \ \ \ \  \sin (\Pi)  = \mathbf{0}, $\\
Resembling the Euler identity, the
  following identity holds:

\begin{equation}\label{BAE2}
 e^{J \ \wedge \ \Pi}  = \mathbf{1} .
\end{equation}

\end{remark}

\begin{remark}(\cite{Memo})
As an application of (\ref{BAE}), we give a shorter proof 
for the fact that:
$  \ \ \cos (X \vee Y) \ \vee \ [ \ J \wedge \sin (X \vee Y) \ ] =
[ \ \cos (X) \ \vee \ [ \ J \wedge \sin (X)  ] \ ] \wedge
[\ \cos (Y) \ \vee \ [ \ J \wedge \sin (Y)  ] \ ] $.

\end{remark}


\begin{remark}
 As explained in the paper \cite{SM}, the above remark could lead to
solutions for the Yang-Baxter equations.
\end{remark}

\begin{remark} The applications of this approach
are diverse. For example, by  Hopf Algebras theory
(the section on
representative coalgebras), we
 are lead to a coalgebra structure:
$ \ \ \Delta (c) = c \otimes c $, $ \ \eps (c) =1 $;
$ \ \Delta (s) = s \otimes c + c \otimes s + s \otimes s $, $ \ 
\eps (s) =0 $.\\
The current Euler formula leads to the existence of a certain coideal.
\end{remark}

\bigskip

\subsection{  Weak B-rings}

Many mathematicians think that
the rings  and the Boolean algebras  are very different structures. 
 However, the study of their axioms (in \cite{Brings})
led to B--rings
   (structures which unify  Boolean algebras and rings).
We will give the definition of a weak B--ring below.

\begin{definition}
 \label{Bring} A { weak \bf B--ring } is a 6--tuple
$ \ ( X, \vee, 0, \cdot , 1 , \hat \ )$, where
$ \ X$ is a set,
$ \vee $ and $\cdot$
are binary operations on $X$,
$  0, 1 \in X$ and
$ \ \hat \ $ is a unary operation on $X$
($ \ \hat x \in X$).\\
These obey the following axioms.

 $ \ ( X, \vee, 0) \ $ is a commutative monoid:\\
(1) $ (x \vee y) \vee z =  x \vee (y \vee z) \ \ \ \ \forall x,y,z \in X$;\\
(2)  $ x \vee y = y \vee x \ \ \ \ \ \ \ \ \ \ \ \forall x,y \in X$;\\
(3) $ x \vee 0 = x \ \ \ \ \ \ \ \ \ \ \ \ \  \ \ \ \ \forall x \in X$.

 $ \ ( X, \cdot , 1 )$ is a monoid:\\
(4) $ (x \cdot y) \cdot z =  x \cdot (y \cdot z)  \ \ \ \ \forall x,y,z \in
X$;\\
(5) $ x \cdot 1 = x = 1 \cdot x \ \  \ \ \forall x \in X  $.

 $0$ is an absorbing element:\\
(6) $ x \cdot 0 = 0 = 0 \cdot x \ \ \ \ \forall x \in X $.

 The second operation is distributive with regard to the 
first operation:\\
(7) $ x \cdot (y \vee z) =  (x \cdot y) \vee (x \cdot z) \ \ \ \forall x,y,z \in
X $;\\
(8) $ (x \vee y) \cdot z =  (x \cdot z) \vee (y \cdot z) \ \ \ \forall x,y,z \in
X $.

 The unary operation has some basic properties:\\
(9) $ \hat{\hat x} = x \ \ \ \ \ \ \ \ \ \ \ \ \ \ \ \ \ \ \ \forall x \in X 
$;\\
(10) $ x \vee \hat x  =  \hat x \vee  x  \ \ \ \ \ \forall x \in X  $;\\
(11) $ x \cdot \hat x  =  \hat x \cdot x \ \ \ \ \ \ \ \  \ \ \ \ \forall x \in
X $.\\

\end{definition}


\begin{remark} 
Let $ \ \mathcal{A} \ = \ ( \ A, \vee, \wedge, \ 0, \ 1, \ \widehat \ ) \ $ 
be a Boolean algebra.

We will consider the following operations on  
$ \mathcal{A}  \times \mathcal{A}  $:

$ \ \ \ \  (a, \ b) {\bf \vee } (c, \ d) = ( a \vee c, \ b \vee d) 
 $, $ \ \ \ \ \widehat{(a, \ b)} = ( \hat{a}, \hat{b}) \ $ and 

$ \ (a, \ b) {\bf \wedge } (c, \ d) = ( (a \wedge c) ,
 \ (a \wedge d) \vee (b \wedge c) \vee (b \wedge d) )  \ $.

Notice that
$ \ \ \ (a, \ b) {\bf \wedge } (1, \ 0) = \ (a, \ b) \ $ and
$ \ \ (0, \ 1) {\bf \wedge } (0, \ 1) = 
( 0, \ 1) $.

We will use the following conventions and notations.
The operation $ \wedge $  has priority over the  operation $ \vee $ , 
and
we will omit the  symbol $ \wedge $
when there is no danger of confusion.
Let $ h = (0, \ 1) $, then
$ \ (a, \ b) = (a,0) \vee (0,b) =
(a,0) \vee [ (0,1) \wedge (b,0)] =
 a \vee h b $.
\end{remark}

\begin{theorem}
 With the above notations,
$ \ \mathcal{A}  \times \mathcal{A} = 
 \ ( \ A \times A , {\bf \vee }, {\bf \wedge }, \ (0,0), \ (1, 0), 
\ \widehat \  \ ) \ $ 
is a weak B-ring.
\end{theorem}

{ \bf Proof. \  }
One could meticulously check all the axioms of a weak B-ring.
The main parts are to verify the associativity of $ \wedge $,
and the distributivity of $ \wedge $ with regard to $ \vee $.
Alternatively, one could observe the analogy with the ring
$ \frac{\mathbb{C}[X]}{X^2 - X} =  \mathbb{C}[h] $ with $ h^2 = h $.
\qed

\bigskip

\subsection{ An Euler formula for weak B--rings}

For a Boolean algebra 
$ \ \mathcal{A} \ = \ ( \ A, \vee, \wedge, \ 0, \ 1, \ \widehat \ ) \ $,
let
 $ \ 0 \neq e \ \in A $.

 We define the ``exponential'' function  
 $ e^{ \ ( \ )} \  : \ A \rightarrow A \ $, by the rule
$ e^x = \widehat{x} \vee e$.

Also, we consider the ``trigonometric'' functions
$ \ C_B \ , S_B \  : \ A \rightarrow A \ $ defined by the rules  
 $ \   C_B (x) = \widehat{x} $ and $ \  S_B (x) = {x} $.

It is easy to check that
\begin{equation}\label{expo}
  e^{x \vee y} = e^{x} \wedge e^{ y} \ ;
\end{equation}
\begin{equation}\label{cos}
 C_B (x \vee y)  =
C_B (x) \wedge C_B ( y) \ ;
\end{equation}
\begin{equation}\label{sin}
 S_B(x \vee y) = [ S_B (x) \wedge C_B ( y) ]  \vee [ C_B (x) \wedge S_B ( y) ] 
\vee [S_B (x) \wedge S_B ( y)] \ ;
\end{equation}

The following formula is a kind of fundamental trigonometric formula:

$$ S_B^{ \ \ 2 \wedge} (X) \vee C_B^{ \ \ 2 \wedge} (X) = 1   \ . $$


Alternatively, the following formula could be considered the 
 fundamental trigonometric formula for this setting:
$$ S_B (X) \vee C_B (X) = 1   \ . $$

\bigskip

The exponential map can be extended on
$ \mathcal{A}  \times \mathcal{A}  $
 as follows.

$ E (a,b) \ = \ ( e^a \wedge \widehat{b}, \ e^a \wedge b ) =
\ [ ( \widehat{a} \vee e) \wedge \widehat{b}, \ ( \widehat{a} \vee e) 
\wedge b ] $.

So, 
$  \ \ E (a \vee h b ) = e^{a \vee h b} = (e^a \ \widehat{b} ) \vee h ( e^a \ b ) =  
e^a \ \widehat{b}  \vee h \ e^a  b  $.

\bigskip

\begin{theorem} (\cite{mMemo})
 With the above notations,

$  \ \ E ((a \vee h b) \ \vee  (x \vee h y)) = 
\ E ((a \vee h b) \ \wedge \ E (x \vee h y) $ .

\end{theorem}


\begin{theorem} (\cite{mMemo})
With the above notations,
 the following Euler type formula holds:
\begin{equation}\label{ABAEd}
 E(h y)  = e^{h  y } =
\ C_B (y) \ \vee \ h  S_B (y) .
\end{equation}

Also, in particular,
\begin{equation}\label{ABAEi}
 e^{  h} = h \ .
\end{equation}
\end{theorem}

\bigskip

\begin{theorem}
 The map $ \phi : \ A \times A \rightarrow \ A \times A $,
defined by the rule 
$$ (x,\ y) \mapsto e^{hy} = ( C_B (y),\ S_B (y)) $$ 

is a solution for
the Braid Condition:
$ \ \ \phi^{12} \ \phi^{23} \ \phi^{12} \ = 
 \  \phi^{23} \ \phi^{12} \ \phi^{23} $.
\end{theorem}

\bigskip
{ \bf Proof. \ \ } The left-hand-side reads:

$  \ \phi^{12} \ \phi^{23} \ \phi^{12} (x, \ y,\ z) =
\ \phi^{12} \ \phi^{23} (\widehat{y}, \ y,\ z) =
\ \phi^{12} \ (\widehat{y}, \ \widehat{z},\ z) =
(z, \ \widehat{z},\ z) $.

The right-hand-side reads:

$  \ \phi^{23} \ \phi^{12} \ \phi^{23} (x, \ y,\ z) =
\ \phi^{23} \ \phi^{12} (x, \  \ \widehat{z},\ z) =
\ \ \phi^{23} \  (z, \  \ \widehat{z},\ z)=
(z, \  \ \widehat{z},\ z) $.

The compassion of the above terms completes the proof.

\qed

\bigskip

Applying the above Euler formula for Boolean maps, leads to the
following theorem.

\begin{theorem} For $ \ x \in A \ $,
 the map $ \phi : \ A \times A \rightarrow \ A \times A $,
defined by the rule 
$$ (a,\ b) \mapsto  ( a \vee xb ,\ ( \hat{x} \vee xa)b \ ) $$ 

is a solution for
the Braid Condition 
($ \ \ \phi^{12} \ \phi^{23} \ \phi^{12} \ = 
 \  \phi^{23} \ \phi^{12} \ \phi^{23} $).
\end{theorem}

\bigskip
{ \bf Proof. \ \ } 
The computations are quite involved. We will only sketch the proof.

The right-hand-side reads:
$  \ \ \ \phi^{23} \ \phi^{12} \ \phi^{23} (a, \ b,\ c) =
\ \phi^{23} \ \phi^{12} (a, \ b \vee x c,\ \hat{x} c \vee x b c) =
\ \phi^{23} \ (a \vee x (b \vee c), 
\hat{x} b \vee x a b \vee x a c , 
\ \hat{x} c \vee x b c) =
( a \vee x (b \vee c), 
\ \hat{x} b \vee x a b \vee x a c \vee x b c , \
\hat{x} c \vee x a b c \ ) $.

The left-hand-side reads:
$  \ \phi^{12} \ \phi^{23} \ \phi^{12} (a, \ b,\ c) =
\ \phi^{12} \ \phi^{23}  (a \vee xb, \ \hat{x}b \vee xab,\ c) =\\
\  =  \ \phi^{12} \   (a \vee xb, 
\hat{x}b \vee xab \vee  x c , 
\ \hat{x} c \vee x a b c \ ) =
( a \vee x (b \vee c), 
\ \hat{x} b \vee x a b \vee x a c \vee x b c , \
\hat{x} c \vee x a b c \ ) $.

\qed

\begin{remark}
 If $ x= 0 $ in the above theorem then 
$ \ \ \phi (a,\ b) =  Id ( a  ,\ b \ ) \ =  ( a  ,\ b \ ).  $
(The identity map is a  solution to the Braid Condition.)

 If $ x= 1 $ in the above theorem then 
$ \ \ \phi (a,\ b) =  P (a,\ b) \  =   ( a \vee b ,\ ab \ ).  $
This is a well-known solution to the Braid Condition (see \cite{nichiF}).

 The above map $ \phi (a,\ b) =  ( a \vee xb ,\ ( \hat{x} \vee xa)b \ ) $,
can be written as 

$$ \ \ \phi  =  
\ C_B (x) \ Id \ \vee \  P \  S_B (x) \ .
$$

Indeed, 
$ \phi (a,\ b) =  
\ C_B (x) (a,b) \ \vee \  ( a \vee b ,\ ab \ ) S_B (x)
$, 
because

$ \ (\hat{x}a, \ \hat{x}b) \ \vee \  [ x(a \vee b) ,\ xab \ ])= 
\ (a \vee xb, \ \hat{x}b \vee \ xab \ ) \ .
$

\end{remark}


\section{Series and   Sequences}

This section is based on discussions at {\em Math Cafe} in 2026.
Started as a high-school problem it led to a sophisticated
approach to unify the arithmetic and geometric series.
The idea is to embed the sum of non-zero digits and a corresponding 
geometric
series into a unique formula. For example, the following matrices
are an illustation for this unification theorem.

\begin{equation} \label{rmatcon2}
A= \begin{pmatrix}
1 & 1 \\
0 & 1 \\
\end{pmatrix}
\end{equation}

\begin{equation} \label{rmatcon2}
B= \begin{pmatrix}
2 & 0 \\
0 & \frac{1}{2} \\
\end{pmatrix}
\end{equation}

\bigskip

\begin{theorem}
 Let $ \ A \in {M}_2 (\mathbb{R}) $ such that $ \ \det (A) = 1 $.

Then, $ \ \sum^9_{k=1} \  A^k = ( tr(A) + 1) ( tr(A^3) + 1) \ A^5 \ $.
\end{theorem}

\section{  APPENDIX. Mathematics and Poetry}

\bigskip

 March is a month which should be dedicated 
to Solomon Marcus (March 1-st, 1925 -- March 17-th, 2016).
I wrote an article in a prestigious literature journal 
about the meetings of Solomon Marcus with Nichita Stanescu.
Below, we give a translation
of a poetry from that article (\cite{F}): 

\bigskip

\begin{center}
 { \em As a noble replica\\
To ``Poetica Mathematica'',\\
The Poet wrote on the appendix\\
''The Poetic Mathematics``.}
\end{center}
 
\bigskip

We present a translation
of the well-known Romanian Spring song. In this
case we are interested in saving the rhythm of the Romanian song,
rather then making a very rigorous translation.

\bigskip

\begin{center}
{ \em  Happy, Happy Celebration\\
Spring is coming from vacation...\\
Let's go out with other kids\\
Let's pick flowers from the fields.}
\end{center}

\end{document}